\documentstyle[12pt]{article}
\begin{document}
\begin{center}
{\bf On the Tangent Space to the Universal Teichm\" uller Space}
\vskip .5em
SUBHASHIS NAG
\end{center}
\begin{abstract}
We find a remarkably simple relationship between the following two models of
the tangent space to the Universal Teichm\" uller Space :\\
(1) The real-analytic model consisting of Zygmund class vector fields on the
unit circle;\\
(2) The complex-analytic model comprising 1-parameter families of schlicht
functions on the exterior of the unit disc which allow quasiconformal
extension.

Indeed, the Fourier coefficients of the vector field in (1) turn out to be
essentially the same as (the first variations of) the corresponding power
series coefficients in (2).

These identities have many applications; in particular, to conformal welding,
to the almost complex structure of Teichm\" uller space,
to study of the Weil-Petersson
metric, to variational
formulas for period matrices, etc.. These utilities are explored. \\
{\sl Subject Classification Numbers:} 32G15, 30F60, 30C60.
\end{abstract}
\baselineskip=24pt

\section{Introduction}

Let $ \Delta $ denote the open unit disc, and $ S^{1} = \partial \Delta $.
Two classic models of the Universal Teichm\" uller Space $T(1) = T(\Delta)$ are
well-known (see [6],[7]) : \\
\noindent (a) the real-analytic model containing all (M\" obius-normalised)
quasisymmetric homeomorphisms of the unit circle $S^1$ ; \\
\noindent (b) the complex-analytic model comprising all (normalized) schlicht
functions on the exterior of the disc :
$$
\Delta^\star = \{ z \in \widehat{{\bf C}}:\;\mid z\mid > 1 \} = \widehat{{\bf
C}} - (\Delta \cup S^1)
$$

\noindent which allow quasiconformal extension to the whole of $\widehat{{\bf
C}}$ (the Riemann Sphere).

\par The relationship between them is via the rather mysterious operation
called ``Conformal Welding" (see [5], and below).  Nevertheless, at the
infinitesimal level, the above models have an amazingly simple relationship
that forms the basis for this paper.  Indeed, the $k^{th}$ {\it{Fourier
Coefficient}} of the vector field representing a tangent vector in model (a),
and the (first variation of) the $k^{th}$ {{\it{power series
coefficient}} representing the {\it{same}} tangent vector in model (b),
turn out to be just ($\sqrt{-1}$ times) complex conjugates of each other. See
Theorem 1 below.

It seems that this unexpectedly simple relationship has several interesting
consequences.
\begin{enumerate}
\item It allows a description of the tangent space to $T(\Delta)$ by ``Zygmund
class power series". (Section 3)
\item It provides immediate proof for the remarkable fact that the Hilbert
transform on Zygmund class vector fields on $S^1$ represents the almost complex
structure on $T(\Delta)$. (Section 4)
\item It provides a simple excplicit formula for the derivative of the
conformal welding map. (Section 5)
\item The infinite-dimensional Weil-Petersson metric on (the ``smooth
points" of) $T(\Delta)$ that was found by the present another in [10] Part II
gets a new expression. \\
(Section 6)
\item We get a formula for the derivative of the infinite-dimensional
period mapping studied by us in [8], [9] in terms of power series variations.
This relates to formulas claimed in [4]. (Section 7)
\end{enumerate}
\center${\star~~~~~\star~~~~~\star}$

\bf Acknowledgement}
\par I would like to thank Clifford Earle and Dennis Sullivan for helpful
communications.  In particular, C.Earle suggested some fairly explicit examples
that led me to my general results.

\section{Teichm\" uller Theory}

\par The universal Teichm\" uller space $ T(\Delta) $ is a holomorphically
homogeneous complex Banach manifold that serves as the universal ambient space
where all the Teichm\" uller spaces (of arbitrary Fuchsian groups) lie
holomorphically embedded.
\par As usual, we set the stage by introducing the chief actor -- namely the
space of (proper) Beltrami coefficients $ L^{\infty}(\Delta)_{1} $ ; it is the
open unit ball in the complex Banach space of $L^{\infty}$ functions on the
unit disc $\Delta$.  The principal construction is to solve the Beltrami
equation.
$$
w_{\bar{z}}  =  \mu w_{z} \eqno(1)
$$

\noindent for any $\mu \in L^{\infty} (\Delta)_{1}$.  The two
above-mentioned models of Teichm\" uller space correspond to discussing two
pertinent solutions for (1) :

\noindent {\bf (a)} $w_{\mu}$ - theory : The quasiconformal hemeomorphism of
{{\bf C}} which is $\mu$-conformal (i.e. solves (1)) in $\Delta$, fixes $\pm 1$
and $-i$, and keeps $\Delta$ and $\Delta^{\star}$ (= exterior of $\Delta$) both
invariant.  This $w_{\mu}$ is obtained by applying the existence and uniqueness
theorem of Ahlfors-Bers (for (1)) to the Beltrami coefficient which is $\mu$ on
$\Delta$ and extended to $\Delta^{\star}$ by reflection $(\tilde{\mu}
(1/\bar{z}) = \overline{\mu(z)} z^{2}/\bar{z}^{2}$ for $z \in \Delta)$.

\noindent {\bf (b)} $w^{\mu}$ - theory : The quasiconformal homeomorphism on
{\bf C}, fixing $0,1,\infty$, which is $\mu$-conformal on $\Delta$ and
conformal on $\Delta^{\star}$. $w^{\mu}$ is obtained by applying the
Ahlfors-Bers theorem to the Beltrami coefficient which is $\mu$ on $\Delta$ and
zero on $\Delta^{\star}$.
\par The fact is that $w_{\mu}$ depends only real analytically on $\mu$,
whereas $w^{\mu}$ depends complex-analytically on $\mu$.  We therefore obtain
two standard models ({\bf{(a)}} and {\bf{(b)}} below) of the universal
Teichm\" uller space, $T(\Delta)$.

Define the universal Teichm\" uller space :
$$
T(\Delta) = {L}^{\infty} (\Delta)_{1} / \sim . \eqno(2)
$$

\noindent Here $\mu \sim \nu$ if and only if $w_{\mu} = w_{\mu}$ on $\partial
\Delta = S^{1}$, and that happens if and only if the conformal mappings
$w^{\mu}$ and $w^{\nu}$ coincide on $\Delta^{\star} \cup S^{1}$.

We let
$$
\Phi : L^{\infty}(\Delta)_{1} \longrightarrow T(\Delta) \eqno(3)
$$

\noindent denote the quotient (``Bers") projection.  $T(\Delta)$ inherits its
canonical structure as a complex Banach manifold from the complex structure of
$L^\infty(\Delta)_{1}$ ; namely, $\Phi$ becomes a {\sl{holomorphic
submersion}}.

The derivative of $\Phi$ at $\mu = 0$ :
$$
d_{0} \Phi : L^\infty(\Delta) \longrightarrow T_{\sl{O}} T(\Delta) \eqno(4)
$$

\noindent is a complex-linear surjection whose kernel is the space $N$ of
``infinitesimally trivial Beltrami coefficients''.
$$
N = \{ \mu \epsilon L^{\infty}(\Delta) : \int \int_{\Delta} \mu \phi = 0~~
\mbox{ for all}~~~\phi \in A(\Delta) \} \eqno(5)
$$

\noindent where $A(\Delta)$ is the Banach space of integrable $(L^1)$
holomorphic functions on the disc.  Thus, the tangent space at the origin
$(\sl{O}= \Phi(0))$ of $T(\Delta)$ is $L^{\infty}(\Delta)/N$.

\par See Ahlfors[2], Lehto[6], and Nag[7] for this material and for what
follows.

\par It is now clear that to $\mu \epsilon L^{\infty}(\Delta)_{1}$ we can
associate the {\it{quasisymmetric homeomorphism}}
$$
f_{\mu} = w_{\mu} \mid_{S^{1}} \eqno(6)
$$

\noindent as representing the Teichm\" uller point $[\mu]$ in version
{\bf{(a)}} of $T(\Delta)$.  Indeed $T(\Delta)_{(a)}$ is the homogeneous space :
$$
\begin{array}{llcl} {\bf{(a)}} & T(\Delta) & = & \mbox{{Homeo$_{q.s.}$}
$(S^{1})$ / {M\"ob} $(S^1)$} \\
& & = &  \{ \mbox{quasisymmetric homeomorphisms of $S^1$ fixing $\pm 1$ and
$-i$} \}
\end{array}
$$

Alternatively, $[\mu]$ is represented by the univalent function
$$
f^\mu = w^\mu \mid_{\Delta^{\star}} \eqno(7)
$$

\noindent on $\Delta^{\star}$, in version {\bf (b)} of $T(\Delta)$.  A more
natural choice of the univalent function representing $[\mu]$ is to use a
different normalisation for the solution $w^{\mu}$ (since we have the freedom
to post-compose by a M\" obius transformation).  In fact, let
$$
W^{\mu} = M^\mu o w^\mu \eqno(8)
$$

\noindent where $M^\mu$ is the unique M\" obius transformation so that the
univalent function (representing $[\mu]$) :
$$
F^\mu = W^\mu \mid_{\Delta^{\star}} \eqno(9)
$$

\noindent has the properties: \\
\noindent (i)  $F^\mu$ has a simple pole of residue 1 at $\infty$ \\
\noindent (ii) $(F^{\mu}(\zeta) - \zeta) \rightarrow 0$ as $\zeta \rightarrow
\infty$.

\noindent Namely, the expansion of $F^\mu$ in $\Delta^{\star}$ is of the form :
$$
F^\mu(\zeta) = \zeta (1 + {{c_{2}} \over {\zeta^{2}}} + {{c_{3}} \over
{\zeta^{3}}} + \ldots) \eqno(10)
$$

\noindent Let us note that the original ($0,1,\infty$ fixing) normalisation
gives an expansion of the form :
$$
f^\mu(\zeta) = \zeta (a + {{b_{1}} \over {\zeta}} + {{b_{2}} \over
{\zeta^{2}}} + {{b_{3}} \over {\zeta^{3}}} \ldots) \eqno(11)
$$

\noindent and the M\" obius transformation $M^\mu$ must be $M^{\mu }(w) = w/a -
b_{1}/a$.  Since $(a,b_{1},b_{2},..)$ depend holomorphically on $\mu$, we see
that $(c_{2},c_{3},...)$ also depend holomorphically on $\mu$.  Thus, our
complex-analytic version $T(\Delta)_{(b)}$ of the universal Teichm\" uller
space is :

$$
\begin{array}{lrcl}
{\bf{(b)}} &  T(\Delta) & = &  \left\{{\rm  Univalent\;
functions \;in \;\Delta^{\star}\; with\; power \;series \;of\; the\; form\;
 (10),}\right.  \\
 & & & \left. {\rm allowing\; quasiconformal\; extension\; to\; the\; whole\;
  plane }\right\}.
\end{array}
$$

\noindent $T(\Delta)_{(b)}$ is simply a\,\,``pre-Schwarzian-derivative''\,\,
version of the Bers embedding of Teichm\" uller space.

\par It is worth remarking here that the criteria that an expansion of the form
(10) represents an univalent function, and that it allows quasiconformal
extension, can be written down solely in terms of the coefficients $c_{k}$,
(using the Grunsky inequalities etc.). See Pommerenke [11].  Thus
$T(\Delta)_{(b)}$ can be thought of as a certain space of sequences
$(c_{2},c_{3},...)$, and its tangent space will be given the concomitant
description below.

\noindent {\underline{Tangent space to the real-analytic model}} :  Since
$T(\Delta)$ is a homogeneous space (see version {\bf (a)}) for which the right
translation (by any fixed quasisymmetric homeomorphism) acts as a biholomorphic
automorphism, it is enough in all that follows to restrict attention to the
tangent space at a single point -- the origin ({\sl{O}} = class of the
identity homeomorphism) -- of $T(\Delta)$.

\par Given any $\mu \in L^\infty (\Delta)$, the tangent vector $d_{0}\Phi(\mu)$
is represented by the real vector field $V[\mu] = \dot{w}[\mu] {{\partial}
\over {\partial z}}$ on the circle that produces the 1-parameter flow $w_{t
\mu}$ of quasisymmetric homeomorphisms:
$$
w_{t \mu}(z) = z + t \dot{w} [\mu](z) + o(t) \eqno(12)
$$

\par The vector field becomes in the $\theta$-coordinate :
$$
V[\mu] = \dot{w}[\mu](z) {{\partial} \over {\partial z}} =
u(e^{i \theta}) {{\partial} \over {\partial \theta}} ,
$$
\noindent where,
$$
u(e^{i \theta}) = {{\dot{w}[\mu](e^{i \theta})} \over {i e^{i \theta}}}~.
\eqno(13)
$$

By our normalization, $u$ vanishes at $1,-1$ and $-i$.

In Gardiner-Sullivan $[3, Section8]$ the precise class of vector fields
arising from such quasisymmetric flows is determined as the Zygmund
$\Lambda$ class. They have delineated the theory on the upper half-plane
$U$; we adapt that result to the disc using the M\" obius transformation
$$
T(z) = {{z-i} \over {z+i}}~~,~~T : U \longrightarrow \Delta~~.
\eqno(14)
$$

\noindent We point out that $(0,1,\infty)$ go to $(-1,-i,1)$ respectively.
Notice that the corresponding identification of the real line to $S^1$ is given
by
$$
x = - cot {{\theta} \over {2}}~~, {\rm or}~~,~~ e^{i \theta} = {{x-i} \over
{x+i}}~~. \eqno(15)
$$

\noindent The continuous vector field~~$u(e^{i \theta}) {{\partial} \over
{\partial \theta}}$~~becomes, on~~{\bf R},~~$F(x) {{\partial} \over {\partial
x}}$~~ with
$$
F(x) = {{1} \over {2}} (x^2 + 1) u{({{x-i} \over {x+i}})}~~. \eqno(16)
$$

\noindent Conversely,
$$
u(e^{i \theta}) = {{2 F(x)} \over {x^2 +1}}~~. \eqno(17)
$$

\noindent Since $u$ vanishes at (-1,-i,1), we see
$$
F(0) = F(1) = 0 \mbox{ and}~~{{F(x)} \over {x^2+1}} \rightarrow 0\;{\rm
as}\;x \rightarrow \infty~. \eqno(18)
$$

\noindent Introduce (following Zygmund [13]),
$$
\begin{array}{lll}
\Lambda ({\bf R}) & = & \{ F : {\bf R}
\rightarrow {\bf R}; F {\rm is\; continuous,\; satisfying\; normalizations\;
(18);}  \\
& & {\rm and}, \mid F(x+t) + F(x-t) -2F(x) \mid \leq B\mid t\mid {\rm for\;
some}\; B, \\
& & {\rm  for\; all}\; x\; {\rm and}\; t\; {\rm real}. \}
\end{array}
\eqno(19)
$$

\noindent $\Lambda({\bf R})$ is a (non-separable) Banach space under the
Zygmund norm -- which is the best constant $B$ for $F$. Namely,
$$
|| F || = {\sup_{x,t}} \left \vert {{F(x+t) + F(x-t) - 2F(x)} \over {t}}
\right \vert \eqno(20)
$$

\noindent In [3] it is shown that $\Lambda ({{\bf R}})$ comprises precisely the
vector fields for quasisymmetric flows on ${{\bf R}}$.  Hence, the tangent
space
to version {\bf (a)} of $T(\Delta)$ becomes :
$$
T_{{\sl O}} T(\Delta)_{{(a)}} = \left\{ \begin{array}{lll}  u(e^{i \theta})
{{\partial} \over {\partial \theta}} : &
(i) &  u : S^1 \rightarrow {{\bf R}}~~\mbox{is continuous}~, \\
& & \mbox{vanishing at} (1,-1,-i) ; \\
& (ii) & F_{u}(x) = {{1} \over {2}} (x^2+1) u{({{x-i} \over {x+i}})}~~\mbox{is
in}~~\Lambda({{\bf R}})
\end{array} \right\}
\eqno(21)$$

\noindent {\bf Remark.} The normalization by M\" obius corresponds to adding
an arbitrary $sl(2,{\bf R})$ vector field, $(ce^{i \theta} + \bar{c}e^{-i
\theta} + b) {{\partial} \over {\partial \theta}}~~,~(c \in {{\bf C}}, b \in
{{\bf R}})$, to $u$.  On the real line this is exactly adding an arbitrary real
quadratic polynomial to $F(x)$.  These operations allow us to enforce the
3-point normalization in each description.

We will say a continuous function $u : S^1 \rightarrow$ {{\bf R}} is in
{\it{the Zygmund class}} $\Lambda(S^1)$ {\it{on the circle}}, if, after adding
the requisite $(ce^{i \theta} + \bar{c}e^{i \theta} + b)$ to normalize $u$, the
function satisfies (21). [Can we find a characterization of $\Lambda(S^1)$ in
terms of the decay properties of the Fourier coefficients?]

\underline{Tangent space to the complex analytic model} :

A tangent vector at {\sl{O}} (the identity mapping) to $T(\Delta)_{(b)}$
corresponds to a 1-parameter family $F_{t}$ of univalent functions (each
allowing quasiconformal extension):
$$
F_{t}(\zeta) = \zeta (1 + {{c_{2}(t)} \over {\zeta^{2}}} +  {{c_{3}(t)}
\over {\zeta^{3}}} + \ldots)~,~~{\rm in}~~\mid \zeta \mid > 1~, \eqno(22)
$$

\noindent with $c_{k}(t) = t \dot{c_{k}}(0) + o(t)~, k = 2,3,4, \ldots$.  The
sequences $\{ \dot{c}_{k}(0)~,~k \geq 2 \}$ arising this way correspond
uniquely to the tangent vectors.

\par Our theorem below will allow us to say precisely which sequences occur
(See Corollary 1).

\section{The Promised Relationship}

\par The principal ingredient in the stew is, of course, the infinitesimal
theory for solutions of the Beltrami equation.

\par For any $\nu \in L^{\infty}({{\bf C}})$ let $w^{t\nu}$ be the
quasiconformal homeomorphism of the plane, fixing $0,1,\infty$, and having
complex dilatation (i.e., Beltrami coefficient) $t \nu$ ; (t small complex).
Then, (see, for example, Ahlfors [2, p.104]), uniformly on compact $\zeta$-sets
we have
$$
\begin{array}{rcl}
w^{t \nu} (\zeta) & = & \zeta + t \dot{f}(\zeta) + o(t)~, (t \rightarrow o) \\
\dot{f}(\zeta) & = & - {{\zeta(\zeta -1)} \over {\pi}} \int\!\int_{\bf C} \quad
{{\nu (z)} \over {z(z-1)(z-\zeta)}} dx dy
\end{array}
\eqno(23)
$$

\par For version {\bf (b)} considerations, apply this to
$$
\nu = \left\{ \begin{array}{ll}
\mu & \mbox{on $\Delta$} \\
0   & \mbox{on $\Delta^{\star}$}
\end{array}\right.
\eqno(24)$$

\noindent We see that
$$
{{\partial} \over {\partial t}} \left\vert_{t=0} {(f^{t \mu} (\zeta))} = -
{{\zeta(\zeta -1)} \over {\pi}} \int\!\int_{\Delta} \quad {{\mu (z)} \over
{z(z-1)(z- \zeta)}} dx dy~,\;\mid \zeta \mid > 1~, \right.
\eqno(25)
$$

\noindent with the univalent functions $f^{t \mu}$ as in (11) above.  Expand
$(z- \zeta)^{-1}$ in powers of $\zeta^{-1}$~, collect terms, and compare with
$$
f^{t \mu} (\zeta) = \zeta \left( a(t) + {{b_{1}(t)} \over {\zeta}} +
{{b_{2}(t)}
\over {\zeta^{2}}} + \ldots \right)~~. \eqno(26)
$$

\noindent One obtains (dot represents ${{\partial} \over {\partial t}}$) :
$$
\begin{array}{lcl}
\dot{a}(0) & = & {{1} \over {\pi}} \int\!\int_{\Delta} {{\mu (z)} \over
{z(z-1)}} dx dy \\
\dot{b}_{k}(0) & = & {{1} \over {\pi}} \int\int_{\Delta}~~{\mu (z) z^{k-2}}
dx dy~~,~~k \geq 1
\end{array}
\eqno(27)$$

\noindent The associated normalized univalent functions
$$
F^{t \mu}(\zeta) = \zeta \left( 1 + {{c_{2}(t)} \over {\zeta^{2}}} +
{{c_{3}(t)} \over {\zeta^{3}}} + \ldots \right)~~, \eqno(28)
$$

\noindent have coefficients $c_{k}(t) = b_{k}(t)/a(t)$.  Consequently, we
derive easily (since $a(0) = 1~, b_{k}(0) = 0$ ) :
$$
\dot{c}_{k}(0) = \dot{b}_{k}(0) = {{1} \over {\pi}} \int\int_{\Delta}~~{\mu
(z) z^{k-2}} dx dy~~,~~k \geq 2 \eqno(29)
$$

\par Our aim is to compare these formulas with the Fourier coefficients of the
vector field $V[\mu]$ corresponding to the same $\mu$ in version {\bf (a)}.
Applying (23) to
$$
\nu = \left\{ \begin{array}{ll}
\mu & \mbox{on $\Delta$} \\
\tilde{\mu} & \mbox{(obtained by ``reflection'' of $\mu$) on $\Delta^{\star}$},
\end{array}
\right. \eqno(30)
$$

\noindent and keeping track of the normalizations, one gets (compare p.134 of
[10, Part II]) :
$$
w_{t \mu}(\zeta) =  \zeta + t \dot{w} [\mu](\zeta) + o(t)~~, t
\rightarrow 0,
$$
$$
\begin{array}{lcll}
\dot{w}[\mu](\zeta) & = & - {{(\zeta -1)(\zeta +1)(\zeta +i)}
\over {\pi}} & \left\{ \int\!\int_{\Delta}\quad {{\mu(z)} \over
{(z-1)(z+1)(z+i)(z-\zeta)}} dxdy \right. \\
&  &  & + \left. i {\int\!\int_{\Delta}} \quad {{\overline{\mu(z)}} \over
{(\bar{z}-1)(\bar{z}+1)(\bar{z}-i)(1-\zeta \bar{z})}} dx\,dy \right\}
\end{array}
\eqno(31)
$$

\noindent Now we want to expand in Fourier series the vector field $V[\mu]$ :
$$
u(e^{i \theta}) = {{\dot{w}[\mu](e^{i \theta})} \over {ie^{i \theta}}} =
\sum_{k=-\infty}^{\infty} a_{k} e^{ik \theta} \eqno(32)
$$

\noindent Since $u$ is real valued, one knows $a_{-k} = \bar{a}_{k}, \quad k
\geq 1$. Calculating the $a_{k}$ from (31) one derives, after taking care of
some remarkable simplifications, (to which I drew attention in [10] also), that
$$
a_{-k} = - {{i} \over {\pi}} \int\!\int_{\Delta} \quad \mu(z) z^{k-2} dx \,
dy,\quad k \geq 2~~.
\eqno(33)
$$

\noindent (The remark after (21) shows that $a_{o}$ and $a_{1}$ do not matter
owing to the $sl(2,{{\bf R}})$ normalization.)

{\bf Theorem 1.}
The tangent vector to $T(\Delta)$ represented by $\mu \in L^\infty (\Delta)$,
corresponds to the Fourier expansion
$$
u(e^{i \theta}) = \sum_{k = -\infty}^{\infty} a_{k} e^{ik\theta} \quad
{\mbox{in version {\bf (a)}}}~~.
$$

\noindent The same $\mu$ corresponds to the 1-parameter family of schlicht
functions
$$
F^{t \mu} (\zeta) = \zeta \left( 1 + {{t \dot{c}_{2}(0)} \over
{\zeta^{2}}} + {{t \dot{c}_{3}(0)} \over {\zeta^{3}}} + \ldots \right) + o(t)
$$

\noindent in version {\bf (b)}.  The identities
$$
{\fbox{ $\dot{c}_{k}(0) = i a_{-k} = i \bar{a}_{k}$, \quad {\mbox{for every}}
\quad {$k \geq 2$}}}
\eqno(34)
$$

\noindent hold.

{\bf Proof :} Compare (29) with (33).

We now have the promised precise description of the complex-analytic tangent
space.

\noindent{\bf Corollary 1} As we saw at the end of Section 2, a tangent vector
to $T(\Delta)_{(b)}$ is determined by a complex sequence $\left(
\dot{c}_{2}(0), \dot{c}_{3}(0), \ldots \right)$.  Precisely those sequences
$(\gamma_{2}, \gamma_{3}, \ldots)$ occur for which the function
$$
u(e^{i \theta}) = i \sum_{k=2}^{\infty} \bar{\gamma}_{k} e^{ik\theta} -i
\sum_{k=2}^{\infty} \gamma_{k} e^{-ik\theta} \eqno(35)
$$

\noindent is in the Zygmund class on $S^1$.

{\bf{Harmonic (Bers') Beltrami Coefficients }}:

Introduce the Banach space of Nehari-bounded holomorphic ``quadratic
differentials''
$$
B(\Delta) = \left\{ \phi \in {Hol}(\Delta) : || \phi (z) (1- \mid z
\mid^2)^2 ||_\infty < \infty \right\} \eqno(36)
$$

\noindent To every $\phi \in B(\Delta)$ we associate the $L^\infty$ function
$$
\nu_{\phi} = \overline{\phi(z)} (1-\mid z \mid^2)^2~~{\rm on}~~\Delta.
\eqno(37)
$$

\noindent Foundational results of Teichm\" uller theory guarantee that the
Beltrami coefficients $\{\nu_{\phi} : \phi \in B(\Delta) \}$ form a
complementary subspace to the kernel $N$ (see equation (5)) of the map
$d_{0}\Phi$.  Thus, this space of {\sl{harmonic Beltrami coefficients}} :
$$
H = \{ \nu \in L^{\infty} (\Delta) : \nu \mbox{ is of the form (37) for some}~~
\phi \in B(\Delta) \} \eqno(38)
$$

\noindent is isomorphic to the tangent space $T_{\sl{O}} T(\Delta)$.  In fact,
this remains true for the Teichm\" uller space $T(G)$ of any Fuchsian group $G$
acting on $\Delta$, simply by replacing $B(\Delta)$ by the subspace
$B(\Delta,G)$ comprising those functions which are quadratic differentials for
$G$.  See Ahlfors [2, Chapter 6] and Nag[7, Chapter 3] for all this.

Therefore the tangent vector associated to an arbitrary $\mu$ is also
represented by a unique Beltrami form of harmonic type (37).  For harmonic
$\mu$ we get a beautifully simple formula for the Fourier coefficients, and
hence using (34) also for the power series coefficients, representing that
tangent vector $d_{0} \Phi(\mu)$.

{\bf Proposition 1.}  Let $\mu = \overline{\phi(z)} (1-\mid z \mid^2)^2$ on
$\Delta$ with $\phi \in B(\Delta)$~,
$$
\phi(z) = h_o + h_1z + h_2z^2 + \ldots, in \mid z \mid < 1~~. \eqno(39)
$$

\noindent The relevant Fourier coefficients $a_{k}$ of the corresponding
Zygmund class vector field $V[\mu]$ on $S^1$ are
$$
\bar{a}_{-k} = a_k = {{2i} \over {(k^3-k)}} h_{k-2}~~, {\rm for}~~k \geq 2~~.
\eqno(40)
$$

{\bf Proof :}  Compute using formula (33) above.

{\bf Remark 1.}  In the presence of a Fuchsian group $G,~\phi$ is a (2,0)-form
for $G$ and the vector field $V[\mu]$ is also $G$-invariant.  That imposes
conditions on the coefficients $h_k$ and $a_k$ respectively, which interact
closely with the relationship (40) exhibited above.

{\bf Remark 2.} The Proposition above can be utilised to analyse why Bers
coordinates are geodesic for the Weil-Petersson metric (Section 6). (Vide
Ahlfors[1] and later work of Royden and Wolpert.)

{\bf{Explicit family of examples}} :

Here is a computable family of examples for which $w^{t \mu}$ can be explicitly
determined, and hence our result can be checked.  These examples are a modified
form of some that were suggested to me by Clifford Earle.

Look at $\mu \in L^\infty (\Delta)$ given by
$$
\mu(z) = -n z^2 \bar{z}^{n-1} \eqno(41)
$$

\noindent Here $n \geq 3$ ; ($n=2$ works also, with minor changes).

The vector field $V[\mu]$ on $S^1$ has Fourier coefficients (from
(33)) as exhibited :
$$
a_{k} = \left \{ \begin{array}{ll}
$-i$ & \mbox{for $k=n-1$} \\
$i$ & \mbox{for $k=1-n$} \\
0 & \mbox{for any other $k \geq 2$ or $\leq -2$.}
\end{array} \right.
\eqno(42)
$$

\noindent The interesting thing is that we can explicitly describe $w^{t\mu}$
for these $\mu$, for all complex $t$ satisfying $\mid t \mid < {1\over n}~$.
Indeed, $w^{t\mu}(\zeta) = f^{t\mu}(\zeta)/(1+t)$, where :
$$
f^{t \mu}(\zeta) = \left\{ \begin{array}{ll}
\zeta \left( 1 + {{t} \over {\zeta^{n-1}}} \right)^{-1} & \mbox{on $\mid
\zeta \mid \geq 1$} \\
\left( {{1} \over {\zeta}} + t \bar{\zeta}^{n} \right)^{-1} & \mbox{on $\mid
\zeta \mid \leq 1$}
\end{array} \right.
\eqno(43)
$$

\noindent It is not hard to check that $f^{t \mu}$ is quasiconformal on {{\bf
C}}, and that its complex dilatation on $\Delta$ is $t \mu$~. The $\{ \mid
\zeta \mid \geq 1 \}$ portion in (43) represents the 1-parameter family of
schlicht functions, and we can write down immediately :
$$
\dot{c}_{k}(0) = \left\{ \begin{array}{ll}
$-1$ & \mbox{for $k=n-1$} \\
0 & \mbox{for any other $k \geq 2$.}
\end{array} \right.
\eqno(44)
$$

\noindent This corroborates Theorem 1.

{\bf Remark.} In constructing these examples, it is, of course, the
quasiconformal homeomorphisms (43) that were written down first; (41) and
(42) were derived from it.

\section{The Almost Complex Structure}

Using Theorem 1 we get an immediate proof of the fascinating fact that the
almost complex structure of $T(\Delta)$ transmutes to the operation of Hilbert
transform on Zygmund class vector fields on $S^1$.  {\sl{Namely, we want to
prove that the vector field}} $V[\mu]$ (equation (13)) {\sl{is related to}}
$V[i \mu]$ {\sl{as a pair of conjugate Fourier series}}.

But the tangent vector represented by $\mu$ in the complex-analytic description
$T(\Delta)_{(b)}$ corresponds to a sequence $(\dot{c}_{2}(0), \dot{c}_{3}(0),
\ldots )$, as explained.  Since the $c_{k}$ are holomorphic in $\mu$, the
tangent vector represented by $i \mu$ corresponds to $(i \dot{c}_{2}(0),
i \dot{c}_{3}(0), \ldots)$.  The relation (34) of Theorem 1 immediately shows
that the $k^{th}$ Fourier coefficient of $V[i \mu]$ is $-i.sgn(k)$ times the
$k^{th}$ Fourier coefficient of $V[\mu]$.  We are through.

{\bf Remark 1.} The Hilbert transform description of the complex structure
on the tangent space of the Teichm\" uller space was first pointed out by
S.Kerckhoff.  A proof of this was important for our previous work, and appeared
in [10, Part I].

{\bf Remark 2.} The result above gives an independent proof of the fact that
conjugation of Fourier series preserves the Zygmund class $\Lambda(S^1)$.  That
was an old theorem of Zygmund [13].

\section{Conformal Welding and its derivative}

The Teichm\" uller point $[\mu]$ in $T(\Delta)_{(b)}$ is the univalent function
$F^{\mu}$ [or, equivalently its image quasidisc $F^{\mu}(\Delta^{\star})$].
The same $[\mu]$ appears in $T(\Delta)_{(a)}$ as the quasisymmetric
homeomorphism $w_{\mu}$ on $S^{1}$.  The relating map is the ``Conformal
Welding"
$$
{\bf W} : T(\Delta)_{(b)} \longrightarrow T(\Delta)_{(a)}~~.
$$

\noindent Namely, given a simply connected Jordan region $D$ on
$\widehat{{\bf{C}}}$ one looks at any Riemann mapping $\rho$ of $\Delta$ onto
$D^{\star}$(=exterior of $D$) and also a Riemann mapping $\sigma$ of
$\Delta^{\star}$ onto $D$.  Both $\rho$ and $\sigma$ extend continuously to the
boundaries to provide two homeomorphisms of $S^{1}$ onto the Jordan curve
$\partial D$.  We define the welding homeomorphism :
$$
{\bf W} (D) = \rho^{-1} o \sigma : S^1 \rightarrow S^1~~. \eqno(45)
$$

\noindent We can normalize by a M\" obius transformation so that {\bf W}$(D)$
fixes $1,-1,-i$.

Clearly, since $\rho = w^{\mu} o w^{-1}_{\mu}$ on $\Delta$, and $\sigma =
w^\mu$ on $\Delta^{\star}$, work as Riemann maps, we see that $T(\Delta)_{(b)}$
and $T(\Delta)_{(a)}$ are indeed related by this fundamental operation of
conformal welding.

{\bf Theorem 2.} The derivative at the origin of $T(\Delta)$ to the map {\bf W}
is the linear isomorphism :
$$
d_{\sl{O}}{\bf W} : \{ \left( \dot{c}_{2}(0), \dot{c}_{3}(0), \ldots \right)
\} \rightarrow {\rm Zygmund\; class}\; \Lambda(S^1) \eqno(46)
$$

\noindent sending $(\gamma_2, \gamma_3, \ldots)$ to the vector field $u(e^{i
\theta}) {{\partial} \over {\partial \theta}}$, where
$$
u(e^{i \theta}) = i \sum_{k=2}^{\infty} \bar{\gamma}_{k} e^{ik\theta} -i
\sum_{k=2}^{\infty} {\gamma}_{k} e^{-ik\theta}
\eqno(47)
$$

{\bf Proof :} Follows from Theorem 1 and the remarks above.

{\bf Remark.} Conformal welding has been studied by many authors even for
domains more general than quasidiscs.  See, for example,
Katznelson-Nag-Sullivan[5].  The derivative formula should now be extended to
the larger context.

\section{Diff($S^1$)/M\" ob($S^1$) inside $T(\Delta)$}

As usual, let Diff$(S^1)$ denote the infinite dimensional Lie group of
orientation preserving $C^\infty$ diffeomorphisms of $S^1$.  The complex-
analytic homogeneous space (see Witten [12])
$$
M = {\rm Diff} (S^1) / {\rm Mob}(S^1)
\eqno(48)
$$

\noindent injects holomorphically into $T(\Delta)_{(a)}$.  This was proved in
[10, Part I].  The submanifold $M$ comprises the ``smooth points" of
$T(\Delta)$ ; in fact, in version (b), the points from $M$ are those quasidiscs
$F^{\mu}(\Delta^{\star})$ whose boundary curves are $C^\infty$.

$M$, together with its modular group translates, foliates $T(\Delta)$ -- and
the fundamental Kirillov-Kostant K\" ahler (sympletic) form (Witten[12]) exists
on each leaf of the foliation.  Up to an overall scaling this homogeneous K\"
ahler metric gives the following pairing $g$ on the tangent space at {\sl{O}}
to $M$ :
$$
g(V,W) = Re \left[ \sum_{k=2}^{\infty} a_{k} \bar{b}_{k} (k^3-k) \right]
\eqno(49)
$$
$$
\begin{array}{lclcl}
V & = & \sum_{2}^{\infty} a_{k} e^{ik \theta} & + & \sum_{2}^{\infty}
\bar{a}_{k} e^{-ik \theta}, \\
W & = & \sum_{2}^{\infty} b_{k} e^{ik \theta} & + & \sum_{2}^{\infty}
\bar{b}_{k} e^{-ik \theta},
\end{array}
$$

\noindent represent two smooth vector fields on $S^1$.

The metric $g$ on $M$ was proved by this author [10, Part II] to be {\sl{the
(regularized) Weil-Petersson metric (WP) of universal Teichm\" uller space}}.
Theorem 1 allows us to express the pairing for $g = WP$ in terms of 1-parameter
flows of schlicht functions.

{\bf Theorem 3.} Let $F_{t}(\zeta)$ and $G_{t}(\zeta)$ denote two curves
through origin in $T(\Delta)_{(b)}$ of the form (22), representing two tangent
vectors, say $\dot{F}$ and $\dot{G}$.  Then the Weil-Petersson pairing assigns
the inner product
$$
WP(\dot{F},\dot{G}) = -Re \left[\sum_{k=2}^{\infty} \overline{\dot{c}_{k}(0)}
\dot{d}_k(0) (k^3-k) \right], \eqno(50)
$$

\noindent where $c_k(t)$ and $d_k(t)$ are the power series coefficients for the
schlicht functions $F_t$ and $G_t$ respectively, in equation (22).  The series
above converges absolutely whenever the corresponding Zygmund class functions
(see equations (35) or (47)) are $C^{3/2 + \epsilon}$.

{\bf Proof:} Combine (34) with (49) and the results of [10, Part II].

\section{Variational Formula for the Period Mapping}

Recently in [8],[9] the author has studied a generalisation of the classical
period mappings to the infinite dimensional context of universal Teichm\" uller
space. Indeed, there is a natural equivariant, holomorphic and K\" ahler-
isometric immersion
$$
\prod : M \longrightarrow D_{\infty} \eqno(51)
$$

\noindent of M = Diff($S^1$)/M\" ob($S^1$) into the infinite dimensional
version, $D_{\infty}$, of the Siegel disc.  $D_{\infty}$ is a complex manifold
comprising certain complex symmetric Hilbert-Schmidt ${(Z_{+}\times Z_{+})}$
matrices.  $\prod$ qualifies as a generalised period matrix map, and its
variation satisfies a Rauch-type formula.  See [9].  An extension of $\prod$ to
{\sl{all}} of $T(\Delta)$ is being worked on by this author using ideas
communicated by Dennis Sullivan.

In the work cited above we proved that for any $\mu$ in $L^\infty(\Delta)$, the
$(rs)^{th}$-entry of the period matrix satisfies $(r,s \geq 1)$
$$
\prod ([t \mu])_{rs} = t \sqrt{-rs}\> a_{-(r+s)} + O(t^2) \eqno(52)
$$

\noindent as $t \rightarrow 0$.  Here $a_k$ as usual denotes the Fourier
coefficients of the vector field represented by $\mu$ (equation (32)).

By Theorem 1 we see that the variational formula above may be written
$$
\prod ([t \mu])_{rs} = \sqrt{rs}\> c_{r+s}(t) + O(t^2) \eqno(53)
$$

\noindent where $c_k(t)$ are the power series coefficients appearing in (22)
for the schlicht functions $F^{t \mu}$.  Equation (53) may be compared with
the formula [(30) in their paper] claimed by Hong-Rajeev[4] in just this
setting.
\center${\star~~~~~\star~~~~~\star}$

\baselineskip = 24pt
\begin{center}
{\bf Appendix}\footnote{Appendix to the paper ``On the Tangent
Space to the Universal Teichm\" uller Space'' by SUBHASHIS NAG.} \\
{\bf{A FUGUE ON AHLFORS [1], AHLFORS-WEILL AND OUR RESULTS}}
\end{center}

The intriguing formulas that Ahlfors exhibited (more than thirty years
ago!) in his paper [1], do have a surprising way of cropping up in
various contexts in later work of many authors. Upon reading the present
article, Prof.Clifford Earle has pointed out to me that the formulas
(1.20) and (1.21) in Section 1.6 of Ahlfors [1] could be used to derive an
elegant proof of our main theorem (equation (34)) above.

The crucial thing is to verify our formulas (29) and (33) (in Section 3
above). As explained with equations (37) - (38), every Beltrami form
$\mu$ is uniquely the sum of an infinitesimally trivial one (equation
(5)) and a harmonic form (equation (5)) and a harmonic form (equation
(37)). Thus the formulas need to be checked only for these two types of
forms.  Of course, all the relevant quantities are zero for
infinitesimally trivial forms. {\it The gist of the matter is that for
harmonic Beltrami forms the Ahlfors-Weill section implies formula}
(29), {\it whereas} (1.21) {\it of Ahlfors} [1] {\it implies formula}
(33)!

{\bf Ahlfors-Weill and formula (29):} \hspace{.5cm} Let
$$
\mu = \overline{\phi(z)} (1-|z|^{2})^{2} \eqno(54)
$$

\noindent be a harmonic Beltrami form on $\Delta$ with $\phi$ as in
(39).  The Ahlfors-Weill theorem tells us explicitly the schlicht
function $w^{t\mu}$ on $\Delta^{\star}$ (for small t) and hence allows
us to compute the variation of the power series coefficients,
$\dot{c}_{k}(0)$.  We refer to Section 3.8 of Nag[7] -- especially
3.8.6 -- for that result.

In fact, let $v_{1}$ and $v_{2}$ be linearly independent solutions in
the unit disc of the differential equation:
$$
v'' = \phi v \eqno(55)
$$

\noindent normalized so that $v_{1}(0) = v'_{2}(0) = 1$ and $v_{2}(0) =
v'_{1}(0) = 0$. Then Ahlfors-Weill tells us that (up to possibly a M\"
obius transformation)
$$
w^{t\mu} (\zeta) = \bar{v}_{1} (1/\bar{\zeta}) / \bar{v}_{2}
(1/\bar{\zeta}) \ {\rm for} \ |\zeta| > 1\,. \eqno(56)
$$

Solving (55) for $v_{1}$ and $v_{2}$ by the method of indeterminate
coefficients yields:
$$
v_{1}(z) = 1 + t \sum_{k=2}^{\infty} {{h_{k-2}} \over {k(k-1)}} z^{k}
\quad + o(t) \eqno(57)
$$
$$
v_{2}(z) = z + t \sum_{k=2}^{\infty} {{h_{k-2}} \over {k(k+1)}} z^{k+1}
+ o(t) \eqno(58)
$$

\noindent Substituting these in (56) we deduce quickly that
$$
\dot{c}_{k}(0) = {{2} \over {(k^{3}-k)}} \overline{h_{k-2}}\,, {\rm for}\;
k \geq 2\,. \eqno(59)
$$

\noindent For the harmonic form $\mu$, formula (59) is exactly formula
(29).

\noindent{\bf Remark:} In 3.8.5 of Nag[7] a new proof of the
Ahlfors-Weill theorem was given using an idea of Royden. The
calculations made there are closely relevant to proving (29) {\it
directly} for harmonic $\mu$ -- without passing to series expansions.

{\bf Ahlfors[1] and formula (33):} Formula (1.21) in Ahlfors [1] in our
notation reads:
$$
\dot{w}[\mu] = \bar{\Phi}'' (1-|z|^{2})^{2} + 2 \bar{\Phi}' z(1-|z|^{2})
+ 2 \bar{\Phi}z^{2} -2 \Phi \eqno(60)
$$

\noindent valid for $|z| \leq 1$, where $\mu$ is the harmonic Beltrami
form (54). Here $\Phi$ (holomorphic in $\Delta$) is related to $\phi$
(of (54)) by
$$
\Phi''' = \phi \eqno(61)
$$

\noindent See Ahlfors[1] formula (1.20) for this.  In order to calculate
the Fourier coefficients $a_{k}$, defined as in (32) above, we only
require (60) on the circle $|z| = 1$. The first two terms of (60)
therefore drop off, and a straightforward computation produces:
$$
a_{k} = {{2i} \over {(k^{3}-k)}} h_{k-2}\,, {\rm for}\; k \geq 2.
\eqno(62)
$$

\noindent But this is exactly formula (40), which is formula (33) for
harmonic Beltrami forms.

Consequently, comparing (59) with (62) proves the complete result.

\begin{flushleft}
\bf{The Institute of Mathematical Sciences} \\
\bf{C.I.T.Campus} \\
\bf{Madras  600 113, INDIA}

\bf{e-mail :} nag@imsc.ernet.in
\end{flushleft}
\end{document}